\documentclass[11pt,onecolumn,a4paper,english]{article}
\usepackage[pdftex]{graphicx}
\usepackage{textpos}
\usepackage{babel,varioref}
\usepackage{fancyhdr}
\usepackage{amsmath}
\usepackage{amsfonts}
\usepackage{array}
\usepackage{color}
\usepackage{colortbl}
\usepackage{amsthm}

   \DeclareGraphicsExtensions{.pdf,.jpeg,.png}

\begin{document}
%
%
\title{The Generalization of the Decomposition of Functions by Energy Operators}
%
%
\maketitle
\author{J.-P.~Montillet \footnote{Dr. J.P. Montillet is a Research Fellow at the Research School of Earth Sciences at the Australian National University}, $j.p.montillet@anu.edu.au$}
\begin{abstract}
\boldmath{This work starts with the introduction of a family of differential energy operators. Energy operators (${\Psi}_{R}^{+}$, ${\Psi}_{R}^{-}$) were defined together with a method to decompose the wave equation in a previous work. Here the energy operators are defined following the order of their derivatives ($\Psi^{-}_k$, $\Psi^{+}_k$, $k=\{0,\pm 1,\pm 2,...\}$). The main part of the work demonstrates for any smooth real-valued function $f$ in the Schwartz space ($\mathbf{S}^{-}(\mathbb{R})$), the successive derivatives of the $n$-th power of $f$ ($n \in \mathbb{Z}$ and $n\neq 0$) can be decomposed using only $\Psi^{+}_k$ (Lemma); or if $f$ in a subset of $\mathbf{S}^{-}(\mathbb{R})$, called $\mathbf{s}^{-}(\mathbb{R})$, $\Psi^{+}_k$ and $\Psi^{-}_k$  ($k\in \mathbb{Z}$) decompose in a unique way the successive derivatives of the $n$-th power of $f$ (Theorem). Some properties of the Kernel and the Image of the energy operators are given along with the development. Finally, the paper ends with the application to the energy function.
}
\end{abstract}


%
%
%
%
%
%
%
%
\section{Introduction}\label{Introduction part}
Two decades ago, an energy operator (${\Psi}_{R}^{-}$) was first defined in \cite{Kaiser90}. Since then, this work has been extensively used in telecommunications (see for example \cite{Bovik93} or \cite{Hamila1999}). The bilinearity properties of this operator were studied in \cite{Boudraa et al.2009}. More recently, the author in \cite{JPMontillet2010} introduced the energy operators (${\Psi}_{R}^{+}$, ${\Psi}_{R}^{-}$) in time and space. This was part of a general method for separating the energy of finite energy functions in time and space with application to the wave equation. Note that calling ${\Psi}_{R}^{+}$ an energy operator is an abuse of language. The term was already used in \cite{JPMontillet2010} as the definition of ${\Psi}_{R}^{+}$ is very close to the Teager-Kaiser energy operator ${\Psi}_{R}^{-}$. 
\newline  This work focuses on the decomposition of a smooth real-valued function $f$ using family of differentiable energy operators based on the energy operators ${\Psi}_{R}^{\pm}$. Throughout this work, $f$ is supposed to be in the Schwartz space $\mathbf{S}^{-}(\mathbb{R})$ defined as:
\begin{equation}
\mathbf{S}^{-}(\mathbb{R}) =\{f \in \mathbf{C}^{\infty}(\mathbb{R}), \qquad {sup}_{t<0} |t^k||\partial_t^j f(t)|<\infty,\qquad \forall k \in \mathbb{Z}^+, \qquad \forall j \in \mathbb{Z}^+ \}
\end{equation}
Note that $\partial_t$ is the derivative related to the variable $t$. $\mathbb{Z}^+$ is the set of positive integers including $0$. In the following, let us call the set $\mathcal{F}(\mathbf{S}^{-}(\mathbb{R}),\mathbf{S}^{-}(\mathbb{R}))$ all functions defined such as $F:$ $\mathbf{S}^{-}(\mathbb{R})$ $\rightarrow$ $\mathbf{S}^{-}(\mathbb{R})$.
%
%
\newline First, two families of differential energy operators $\Psi^{+}_k$ and $\Psi^{-}_k$ ($k=\{0,\pm 1,\pm 2,...\}$) are introduced with the notations following \cite{Maragos1995} and \cite{JPMontillet2010} where $k$ is the degree of their derivatives. Note that for $f$ in $\mathbf{S}^{-}(\mathbb{R})$, $\Psi^{+}_k(f)$ and $\Psi^{-}_k(f)$ ($k=\{0,\pm 1,\pm 2,...\}$) are in $\mathbf{S}^{-}(\mathbb{R})$ ($(\Psi^{\pm}_k)_{k \in \mathbb{Z} }$ $\subseteq$ $\mathcal{F}(\mathbf{S}^{-}(\mathbb{R}),\mathbf{S}^{-}(\mathbb{R}))$. If not explicitly written, any families of operator in $\mathcal{F}(\mathbf{S}^{-}(\mathbb{R}),\mathbf{S}^{-}(\mathbb{R}))$ in this work follow the derivative chain rule property:
%
%
\begin{equation}\label{derivative chain rules0111}
\forall f \in \mathbf{S}^{-}(\mathbb{R}), \qquad \partial_t \Psi_k(f) = \Psi_{k+1}(f) + \Psi_{k-1}(\partial_t f)
\end{equation} 
In addition, the term \emph{decompose} is defined as:
\vspace{1.5em}
\newline $\bold{Definition}$ $1$: for all $f$ in $\mathbf{S}^{-}(\mathbb{R})$, for all $v\in\mathbb{Z}^+-\{0\}$, for all  $n\in\mathbb{Z}^+$ and $n>1$, the family of operators $(\Psi_k)_{k \in \mathbb{Z}}$ (with $(\Psi_k)_{k \in \mathbb{Z}}$ $\subseteq$ $\mathcal{F}(\mathbf{S}^{-}(\mathbb{R}),\mathbf{S}^{-}(\mathbb{R}))$) decomposes $\partial_t^v$$f^n$ in $\mathbb{R}$, if it exists $(N_j)_{j\in \mathbb{Z}^+ \cup \{0\}}$ $\subseteq$ $\mathbb{Z^+}$,  $(C_i)_{i=-N_j}^{N_j}$ $\subseteq$ $\mathbb{R}$, and it exists $(\alpha_j)$ and $l$ in $\mathbb{Z^+}\cup\{0\}$ (with $l<v$) 
\\ such as $\partial_t^v$$f^n = \sum_{j=0}^{v-1} \big(_{j}^{v-1} \big) \partial_t^{v-1-j} f^{n-l} \sum_{k=-N_j}^{N_j} C_k \Psi_k(\partial_t^{\alpha_k}f)$.
\vspace{1.5em}
\\  Definition $1$ is based on the general Leibniz rule for the $n-th$ derivative of a product of functions \cite{BruceWest}. One can define the image $Im(\Psi^{+}_k)$ and kernel $Ker(\Psi^{+}_k)$ (for $k$ in $\mathbb{Z}$) of an energy operator such as:
\begin{equation}\label{Impsik+}
Im(\Psi^{+}_k) = \{\Psi^{+}_k(f) \in \mathbb{R} | \qquad f \in \mathbf{S}^{-}(\mathbb{R}) \}
\end{equation}
%
and
\begin{equation}
Ker(\Psi^{+}_k) = \{f\in \mathbf{S}^{-}(\mathbb{R})|  \qquad \Psi^{+}_k(f)=0 \}
\end{equation}
Obviously, the null function ($f: \mathbb{R}\rightarrow 0$) belongs to $Ker(\Psi^{+}_k)$.
One can define also $Im(\Psi^{-}_k)$ and $Ker(\Psi^{-}_k)$ associated with the family of DEOs $\Psi^{-}_k$ ($k$ in $\mathbb{Z}$).
%
Now, let us define a subset $\mathbf{s}^{-}(\mathbb{R})$ $\subseteq$ $\mathbf{S}^{-}(\mathbb{R})$ such as:
\begin{eqnarray}
\mathbf{s}^{-}(\mathbb{R}) &=& \{f\in \mathbf{S}^{-}(\mathbb{R})|  \hspace{0.2em}\forall\hspace{0.2em} k \in \mathbb{Z}, \hspace{0.2em}  \Psi^{+}_k(f) \neq \{0\} \hspace{0.2em}  \nonumber \\
                   & & \hspace{0.2em} \forall\hspace{0.2em} k \in \mathbb{Z}-\{1\}, \hspace{0.2em}\Psi^{-}_k(f) \neq \{0\} \} \nonumber  
\end{eqnarray}
Note that a possible way to define $\mathbf{s}^{-}(\mathbb{R})$ is:
\begin{equation}
\mathbf{s}^{-}(\mathbb{R}) = \{ f \in  \mathbf{S}^{-}(\mathbb{R}) | f \notin  (\cup_{k \in \mathbb{Z}} Ker(\Psi^{+}_k))\cup(\cup_{k \in \mathbb{Z}-\{1\}} Ker(\Psi^{-}_k))\}
\end{equation}
%
The definition of the subset $\mathbf{s}^{-}(\mathbb{R})$ excludes $\Psi^{-}_1$ as by definition of  this operator $Im(\Psi^{-}_1)$ equal $\{0\}$ for all $f$ in  $\mathbf{S}^{-}(\mathbb{R})$.
Following Definition $1$, the \emph{uniqueness} of the decomposition  in $\mathbf{s}^{-}(\mathbb{R})$ with the families of differential operators can be stated as:
\vspace{1.5em}
\newline $\bold{Definition}$ $2$: for all $f$ in $\mathbf{s}^{-}(\mathbb{R})$, for all $v\in\mathbb{Z}^+-\{0\}$, for all  $n\in\mathbb{Z}^+$ and $n>1$, the families of operators $(\Psi^{+}_k)_{k \in \mathbb{Z}}$ and $(\Psi^{-}_k)_{k \in \mathbb{Z}}$ ($(\Psi^{+}_k)_{k \in \mathbb{Z}}$ and $(\Psi^{-}_k)_{k \in \mathbb{Z}}$$\subseteq$ $\mathcal{F}(\mathbf{s}^{-}(\mathbb{R}),\mathbf{S}^{-}(\mathbb{R}))$) decompose uniquely $\partial_t^v$ $f^n$ in $\mathbb{R}$, if for any family of operators $(S_k)_{k \in \mathbb{Z}}$ $\subseteq$ $\mathcal{F}(\mathbf{S}^{-}(\mathbb{R}),\mathbf{S}^{-}(\mathbb{R})$) decomposing  $\partial_t^v$$f^n$ in $\mathbb{R}$, there exists a unique couple $(\beta_1,\beta_2)$ in $\mathbb{R}^2$ such as: 
\begin{equation}
S_k(f) = \beta_1 \Psi^{+}_k(f) + \beta_2 \Psi^{-}_k(f), \qquad \forall k\in\mathbb{Z} 
\end{equation}
\vspace{1.5em}
\\ The main goal of this work is to give a proof of the following lemma and theorem:
\vspace{1.5em}
\\$\bold{Lemma}$: for $f$ in $\mathbf{S}^{-}(\mathbb{R})$, the family of DEO ${\Psi}_{k}^{+}$ ($k=\{0,\pm 1,\pm 2,...\}$) decomposes the successive derivatives of the $n$-th power of $f$ for $n\in\mathbb{Z}^+$ and $n>1$. 
\vspace{1.5em}
\\$\bold{Theorem}$: for $f$ in $\mathbf{s}^{-}(\mathbb{R})$, the families of DEO ${\Psi}_{k}^{+}$ and ${\Psi}_{k}^{-}$ ($k=\{0,\pm 1,\pm 2,...\}$) decompose uniquely the successive derivatives of the $n$-th power of $f$ for $n\in\mathbb{Z}^+$ and $n>1$. 
%
%
\vspace{1.5em}
\\ The proofs of the lemma and theorem are given for the $n$-th power of $f$ with $n\in\mathbb{Z}^+$ and $n>1$. A discussion takes place for the special case $n=1$ and $n<0$. In addition, we assume a function $f$ of a variable $t$ with values in $\mathbb{R}$, which can be time or one of the dimension in space ($x,y,z$). 
%
%
%
Note that in Section \ref{SectionProperties}, the study of the properties of the images helps to simplify the formulas shown in Lemma and  Theorem.
%
%
%
\\ Finally,  the last part is dedicated to applying the development in the previous sections to the energy function $\mathcal{E}$ defined as:
\begin{equation}\label{EnergyfunDefine02}
\mathcal{E}(f_1(\tau)) = \int_a^{\tau} f_1(t)^2dt < \infty
\end{equation}
with $a$ and $\tau$ in $\mathbb{R}$. $f_1$ is assumed to be in $\mathbf{S}^{-}(\mathbb{R})$ and analytic.
%
%
%
%
\section{Family of Energy Operators}
%
%
%
Following the general description in \cite{Kaiser90}, the general formula of the operator can be written as $P_{-}$ a bilinear form of $\mathbb{R}$ $ \times $ $\mathbb{R}$ $\rightarrow$ $\mathbb{R}$ defined in the real domain for all functions $f$ and $g$ in $C^{\infty}(\mathbb{R})$ as:
\begin{equation}
{P}_{-}[f(t),g(t)]=\frac{1}{2}[ f(t)\partial_t g(t)+g(t)\partial_t f(t)]-\frac{1}{2}[f(t)\partial_t^2g(t)+\partial_t^2f(t)g(t)]
\label{Equation Energy Operator in R}
\end{equation}
Some years after, the authors in \cite{Maragos1995} introduced the $k$-th differential energy operator (DEO):  
\begin{equation}\label{familyDEOdefinition}
{\Psi}_{k}^{-}(f(t)) = \dot{f}(t){f}^{(k-1)}(t)- f(t) {f}^{(k)}(t), \qquad k \in \mathbb{Z}
\end{equation}
Note that $k$ is the order of the operators. Here the derivation is following the variable $t$ and ${f}^{(k)}$ means $\partial_t^k f$. One can see that ${\Psi}_{k}^{-}$ is the quadratic form of the bilinear form $P_{-}$.
%
To explicitly define ${f}^{(k)}$ for all k in $\mathbb{Z}-\{0\}$:
\begin{eqnarray}
f^{(k)}(t) &=& \partial_t^k f(t), \qquad \forall k \in \mathbb{Z}^+ -\{0\} \nonumber \\
f^{(k)}(t) &=& \int_{-\infty}^{t}(\hdots(\int_{-\infty}^{\tau_1}f(\tau_1)d\tau_1)...)d\tau_k, \qquad \forall k \in \mathbb{Z}^- -\{0\}\nonumber \\
f^{(k)}(t) &=& f(t), \qquad  k = 0
\end{eqnarray}
With this definition it is important to underline that we are interested in the function such that:
\begin{eqnarray}
\partial_t(\int_{-\infty}^{t}f(s)ds) &=& f(t) \nonumber \\
\int_{-\infty}^{t}\partial_{t}f(s)ds &=& f(t)
\end{eqnarray}
This explains why we choose $f$ in the Schwartz space $\mathbf{S}^{-}(\mathbb{R})$.
%
%
%
%
Based on the definition given in \cite{JPMontillet2010}, one can define the DEO family in the same way:
\begin{equation}\label{Psik+defdefef}
{\Psi}_{k}^{+}(f) = \dot{f}{f}^{(k-1)} + f {f}^{(k)}, \qquad \forall k \in \mathbb{Z}
\end{equation}
${\Psi}_{k}^{+}$ is also a quadratic form with the same properties as ${\Psi}_{k}^{-}$. It is easy to show that the derivative chain rule in Equation \eqref{derivative chain rules0111} remains the same for the DEO family ${\Psi}_{k}^{+}$.
%
%
\section{Proof of the Lemma and Theorem}\label{LemmasProofs}
$\bold{Lemma}$: for $f$ in $\mathbf{S}^{-}(\mathbb{R})$, the family of DEO ${\Psi}_{k}^{+}$ ($k=\{0,\pm 1,\pm 2,...\}$) decomposes the successive derivatives of the $n$-th power of $f$ for $n\in\mathbb{Z}^+$ and $n>1$. 
\begin{proof}
The general proof is structured via an induction for the existence of the decomposition. First, the study of the four first derivatives of some selected powers of $f$ ($n=\{2, 3\}$) shows how the differential energy operators (or DEOs) are defined for each selected power of $f$. Furthermore, it gives a the method to find the energy operator family to decompose the successive derivatives. The case $n=p$ ($p>1$) ends the induction proof.
The non-uniqueness of the decomposition is justified in a separate paragraph with a counter example.
\vspace{1.5em}
\newline A- $\bold{Existence}$ $\bold{of}$ $\bold{the}$ $\bold{Decomposition}$ $\bold{of}$ $\bold{the}$ $\bold{Power}$ $\bold{of}$ $f$ $\bold{via}$ $\bold{DEOs}$
\vspace{1.5em}
\\ $\bold{Case}$ $n=2$:
\vspace{1.5em}
\\Let $f$ be a function in $\mathbf{S}^{-}(\mathbb{R})$. 
%
%
One can write the coefficient $\partial_t^v f^2$ for $v=\{1,2,3,4\}$ as:
\begin{eqnarray}\label{Coeff2n}
\partial_t f^2 &=& 2 f \partial_t f \nonumber \\
\partial_t f^2 &=& {\Psi}_{1}^{+}(f) \nonumber \\
& & \nonumber \\
\partial_t^2 f^2 &=& 2 (\partial_t f)^2 + 2 f \partial_t^2 f \nonumber \\
\partial_t^2 f^2 &=& \partial_t {\Psi}_{1}^{+}(f) \nonumber \\
\partial_t^2 f^2 &=& {\Psi}_{2}^{+}(f) + {\Psi}_{0}^{+}(\partial_t f) \nonumber \\
& & \nonumber \\
\partial_t^3 f^2 &=& {\Psi}_{3}^{+}(f) +2 {\Psi}_{1}^{+}(\partial_t f) + {\Psi}_{-1}^{+}(\partial_t^2 f) \nonumber \\
& & \nonumber \\
\partial_t^4 f^2 &=& {\Psi}_{4}^{+}(f) +3 {\Psi}_{2}^{+}(\partial_t f) +3 {\Psi}_{0}^{+}(\partial_t^2 f) + {\Psi}_{-2}^{+}(\partial_t^3 f) \nonumber \\
\end{eqnarray}
%
In this example, the successive coefficients are calculated using the derivative chain rule (e.g. Equation \eqref{derivative chain rules0111}). The scalars for each coefficient follow the Pascal's triangle rule and the order of the derivatives. Moreover, it is possible to develop the same type of Pascal's triangle rule to predict not only the scalar coefficients but also the order of the energy operators ($k$) involved for each derivative. 
\begin{figure}[!hbp]
\begin{center}
\includegraphics[width=5in]{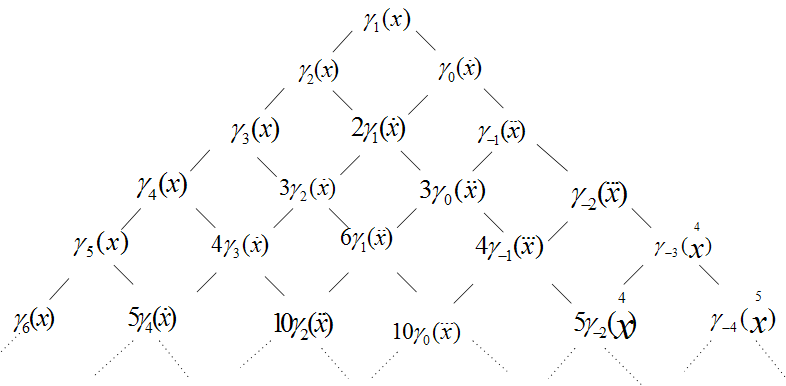}
\caption{\footnotesize{Pascal's triangle rule for the energy operators}}\label{Pascaletrainglerules}
\end{center}
\end{figure}
Figure \ref{Pascaletrainglerules} summarizes this general rule. It is straightforward to see that for the $v$-th derivative order, the highest DEO order is equal to $v$. Then, the remaining energy operator orders are calculated  by decreasing the order by one for each previous  energy operator involved in the ($v-1$)-th derivative. Secondly, it is important to see that all the DEOs involved to approximate the $v$-th derivative of $f^2$ are applied to the $s$-th derivatives of $f$ with $s=\{1,...,v-1\}$, except the DEO with the highest order (which is applied directly to $f$).  
For simplicity using Equation \eqref{Coeff2n}, one can write ($\forall s \in \{1,2,...,m\}$):
\begin{eqnarray}\label{Aocorf1N2}
 a_s^+(f) &=& \partial_t^s f^2 \nonumber \\
a_s^+(f) &=& \partial_t^{s-1} \Psi_{1}^{+} (f) 
\end{eqnarray}
$a_s^+$ is a sum of DEOs. The generalization of $a_s^+(f)$ is given via the formula:
%
%
%
%
%
%
\begin{eqnarray}\label{coefficientsAPgeneral}
a_s^+(f) &=& \sum_{k=0}^{s-1} \big(_{k}^{s-1} \big) \Psi_{2(k+1)-s}^{+}(\partial_t^{s-k-1}f), \qquad \forall s\in \mathbb{Z}^{+}-\{0\} 
\end{eqnarray}
This formula can be checked for $s=\{1,2,3,4\}$ using Equation \eqref{Coeff2n}.
Let us assume the formula is true for $s=m$. One can write:
\begin{eqnarray}\label{coefficientsAPgeneralbis}
a_{m}^+(f) &=& \partial_t a_{m-1}^+(f)\nonumber \\
a_{m}^+(f) &=& \partial_t \sum_{k=0}^{m-2} \big(_{k}^{m-2} \big) \Psi_{2(k+1)-m+1}^{+}(\partial_t^{m-k-2}f)\nonumber \\
a_{m}^+(f) &=& \sum_{k=0}^{m-2} \big(_{k}^{m-2} \big) [ \Psi_{2(k+1)-m+2}^{+}(\partial_t^{m-k-2}f) + \Psi_{2(k+1)-m}^{+}(\partial_t^{m-k-1}f)]
\nonumber \\
\end{eqnarray}
Now using Equations \eqref{coefficientsAPgeneralbis} and the derivative chain rule property, the case $m+1$ is:
%
%
\begin{eqnarray}\label{coefficientsAPgeneralbis2}
a_{m+1}^+(f) &=& \partial_t a_{m}^+(f)\nonumber \\
a_{m+1}^+(f) &=& \sum_{k=0}^{m-1} \big(_{k}^{m-1} \big) [ \Psi_{2(k+1)-m+1}^{+}(\partial_t^{m-k-1}f) + \Psi_{2(k+1)-m-1}^{+}(\partial_t^{m-k}f)] \nonumber \\
a_{m+1}^+(f) &=& \sum_{k=0}^{m} \big(_{k}^{m} \big) \Psi_{2(k+1)-m-1}^{+}(\partial_t^{m-k}f)\nonumber \\
\end{eqnarray}
The above equations finish the induction proof for this particular case.
%
%
%
%
%
\vspace{1.5em}
 \newline $\bold{Case}$ $n=3$:
\vspace{1.5em}
\\ In the same way as in the previous case, we write the development of the first coefficient:
\begin{eqnarray}\label{f3coef1}
\partial_t f^3 &=& 3 f^2 \partial_t f \nonumber \\
\end{eqnarray}
Let us introduce a DEO which is just a product of ${\Psi}_{k}^{+}$ by a constant as:
\begin{equation}\label{GammaOperatorfamily}
\Gamma_{k}^{+} (f) = \frac{3}{2} (\dot{f}{f}^{(k-1)} + f {f}^{(k)}), \qquad \forall k\in \mathbb{Z}
\end{equation}
Clearly, by definition $\Gamma_{k}^{+}$ is a quadratic form and a DEO as it is proportional to ${\Psi}_{k}^{+}$. Note that the derivative properties shown in Equation \eqref{derivative chain rules0111} hold for this DEO. Using Equation \eqref{Aocorf1N2}, it is possible to write a similar equality with $a_k^{+}:$  $A_i^+(f) = 3/2 \partial_t^{i-1} a_i^+(f)$ (with $i$ in $\mathbb{Z}^{+}-\{0\}$). It is then possible to write the Equation \eqref{f3coef1} and the successive derivatives of $f^3$ as :
\begin{eqnarray}
\partial_t f^3 &=& f \Gamma_{1}^{+}(f) \\
\partial_t f^3 &=& f  A_1^{+}(f) \nonumber \\
 & & \nonumber \\
\partial^2_t f^3 &=& f (\Gamma_{2}^{+}(f) + \Gamma_{0}^{+}(\partial_t f)) + \partial_t f \Gamma_{1}^{+} (f) \nonumber \\
\partial^2_t f^3 &=& f A_2^{+}(f)  + \partial_t f A_1^{+}(f) \\
 & & \nonumber \\
\partial^3_t f^3 &=& f (\Gamma_{3}^{+}(f) + 2 \Gamma_{1}^{+}(\partial_t f) + \Gamma_{-1}^{+}(\partial^2_t f)) + 2 \partial_t f (\Gamma_{2}^{+}(f) + \Gamma_{0}^{+}(\partial_t f)) + \partial^2_t f \Gamma_{1}^{+} (f) \nonumber \\
\partial^3_t f^3 &=& f A_3^{+}(f) + 2 \partial_t f A_2^{+}(f) + \partial^2_t f A_1^{+}(f) \\
& & \nonumber \\
\partial^4_t f^3 &=& f (\Gamma_{4}^{+}(f) + 3 \Gamma_{2}^{+}(\partial_t f) + 3 \Gamma_{0}^{+}(\partial^2_t f) + \Gamma_{-2}^{+}(\partial^2_t f)) \nonumber \\
& & + 3 \partial_t f (\Gamma_{3}^{+}(f) + 2 \Gamma_{1}^{+}(\partial_t f) + \nonumber \\
 && \Gamma_{-1}^{+}(\partial^2_t f)) + 3 \partial^2_t f (\Gamma_{2}^{+}(f) + \nonumber \\ 
 & & \Gamma_{0}^{+}(\partial_t f)) + \partial^3_t f \Gamma_{1}^{+} (f) \nonumber \\
\partial^4_t f^3 &=& f A_4^{+}(f) + 3 \partial_t f A_3^{+}(f) + 3 \partial^2_t f A_2^{+}(f) + \partial^3_t f A_1^{+}(f) \nonumber \\
\end{eqnarray}
%
%
\\ There is a certain symmetry between the above equations and Equations \eqref{Coeff2n}. The decomposition of $\partial_t^s f^3$ is performed using the DEO $\Gamma_{k}^{+}$. Using the equations above and the general definition of $A_s^{+}$ ($\forall s \in \mathbb{Z}^+ -\{0\}$) in Equation \eqref{coefficientsAPgeneral}, it is then possible to generalize the formula: 
\begin{equation}\label{dtp3fwithAp}
\partial_t^{m+1} f^3 =\sum_{k=0}^m \big(_{k}^m \big) A_{k+1}^+(f) \partial_t^{m-k} f, \qquad \forall m\in \mathbb{Z}^+
\end{equation}
This formula is checked with $m=\{0,1,2,3\}$. In order to finish the proof by induction, let us then assume the formula true for $m+1$. Following Equation \eqref{dtp3fwithAp} and the previous development, one can write:
\begin{eqnarray}\label{dtp3fwithAp1b}
\partial_t^{m+1} f^3 &=& \partial_t (\partial_t^{m} f^3)\nonumber \\
\partial_t^{m+1} f^3 &=&\sum_{k=0}^{m-1} \big(_{k}^{m-1} \big) [A_{k+2}^+(f) \partial_t^{m-k-1} f +A_{k+1}^+(f) \partial_t^{m-k} f]
\end{eqnarray}
Finally using Equation \eqref{dtp3fwithAp1b}, the development for  $m+2$ is:
\begin{eqnarray}\label{dtp3fwithAp02}
\partial_t^{m+2} f^3 &=& \partial_t(\partial_t^{m+1} f^3)\nonumber \\
\partial_t^{m+2} f^3 &=& \partial_t (\sum_{k=0}^m \big(_{k}^m \big) A_{k+1}^+(f) \partial_t^{m-k} f) \nonumber \\
\partial_t^{m+2} f^3 &=&  \sum_{k=0}^m \big(_{k}^m \big) [A_{k+2}^+(f) \partial_t^{m-k} f +A_{k+1}^+(f) \partial_t^{m-k+1} f] \nonumber \\
\partial_t^{m+2} f^3 &=& \sum_{k=0}^{m+1} \big(_{k}^{m+1} \big) A_{k+1}^+(f) \partial_t^{m-k+1} f \nonumber \\
\end{eqnarray}
This ends the induction proof to show Equation \eqref{dtp3fwithAp}, and also this particular case.
%
%
%
\vspace{1.5em}
\\ $\bold{Case}$ $n=p$ with $p>1$:
\vspace{1.5em}
\\ First, let us assume that the family of DEO ${\theta}_{k}^{+}$ ($k\in\mathbb{Z}$) proportional to the family of DEO  ${\Psi}_{k}^{+}$ ($k\in\mathbb{Z}$) which decomposes the successive derivatives of ${f}^{p-1}$ with the definition:
\begin{eqnarray}
{\theta}_{k}^{+}(f) &=& \frac{(p-1)}{2} (\dot{f}{f}^{(k-1)} + f {f}^{(k)})
\end{eqnarray}
%
In addition, one can define  $B_i^+(f) = \partial_t^{i-1} {\theta}_1^+(f)$ (with $i$ in $\mathbb{Z}^{+}-\{0\}$). Let us then write the first derivatives of ${f}^{p}$:
\begin{eqnarray}\label{NONOVO}
 \partial_t f^p &=&  p f^{p-1} \partial_t f  \nonumber \\
 \partial_t f^p &=&  \frac{p}{2} f^{p-2} \Psi_1^+(f)  \nonumber \\
\partial_t f^p &=& \frac{p}{p-1} {\theta}_{1}^{+}(f) f^{p-2} \nonumber \\
\partial_t f^p &=& \frac{p}{p-1} B_1^+(f)  f^{p-2}\nonumber \\
& & \nonumber \\
%
\partial_t^2 f^p &=& \frac{p}{p-1} B_1^+(f) \partial_t f^{p-2}+ \frac{p}{p-1} B_2^+(f)  f^{p-2} \nonumber \\
\partial_t^3 f^p &=& \frac{p}{p-1}  B_1^+(f) \partial_t^2 f^{p-2}+ 2\frac{p}{p-1}B_2^+(f) \partial_t f^{p-2} \nonumber \\
& & + \frac{p}{p-1}B_3^+(f) f^{p-2}
\end{eqnarray}
As shown in the previous case, we can generalize the formula for the $s$-th derivative as:
\begin{equation}\label{dtppfwithAps}
\partial_t^{s+1} f^p =\sum_{k=0}^s \big(_{k}^s \big) \frac{p}{p-1}B_{k+1}^+(f) \partial_t^{s-k} f^{p-2}, \qquad \forall s \in \mathbb{Z}^+
\end{equation}
The formula in Equation \eqref{dtppfwithAps} has just been verified for $s=\{0,1,2\}$. Furthermore, one can write for the case $s=m$:
\begin{eqnarray}\label{dtp3fwithAp1b}
\partial_t^{m+1} f^p &=& \partial_t (\partial_t^{m} f^p) \nonumber  \\
\partial_t^{m+1} f^p &=&\frac{p}{p-1} \sum_{k=0}^{m-1} \big(_{k}^{m-1} \big) [B_{k+2}^+(f) \partial_t^{m-k-1} f^{p-2} + B_{k+1}^+(f) \partial_t^{m-k} f^{p-2}] \nonumber  \\
\end{eqnarray}
Let us assume that the Equation \eqref{dtppfwithAps} is true for $s=m+1$. Using the previous equation, the case $s=m+2$ is developed as:
\begin{eqnarray}\label{dtp3fwithAp02}
\partial_t^{m+2} f^p &=& \frac{p}{p-1} \partial_t (\sum_{k=0}^m \big(_{k}^m \big) B_{k+1}^+(f) \partial_t^{m-k} f^{p-2}) \nonumber \\
\partial_t^{m+2} f^p &=&  \frac{p}{p-1} \sum_{k=0}^v \big(_{k}^m \big) [B_{k+2}^+(f) \partial_t^{m-k}  f^{p-2} +B_{k+1}^+(f) \partial_t^{m-k+1} f^{p-2}] \nonumber \\
\partial_t^{m+2} f^p &=& \frac{p}{p-1} \sum_{k=0}^{m+1} \big(_{k}^{m+1} \big) B_{k+1}^+(f) \partial_t^{m-k+1} f^{p-2} \nonumber \\
\end{eqnarray}
$(B_k^+)_{k\in  \mathbb{Z}}$ and by definition, $(\theta_k^+)_{k\in  \mathbb{Z}}$ decompose $\partial_t^{s} f^p$  ($s \in \mathbb{Z}^+$, $p\in \mathbb{Z}^+ -\{0,1\}$). This ends the induction proof to confirm Equation \eqref{dtppfwithAps}. 
%
%
%
%
%
Finally as ${\theta}_{k}^{+}$ is proportional to ${\Psi}_{k}^{+}$, one can conclude that the family of DEO  ${\Psi}_{k}^{+}$ ($k\in \mathbb{Z}$) decomposes the successive derivatives of ${f}^{p}$ ($p \in \mathbb{Z}^+$, and $p> 1$). This ends the induction proof on the $n$-th power of $f$. 
%
%
\vspace{1.5em}
\newline B -$\bold{Uniqueness}$ $\bold{of}$ $\bold{the}$ $\bold{Decomposition}$
\\ In the previous proof, it has been shown that given a family of DEOs  ${\Psi}_{k}^{+}$ ($k\in \mathbb{Z}$), it is possible to decompose the successive derivatives $\partial_t^i f^n$ ($n\in \mathbb{Z}^+$, $n>1$, $i\in\mathbb{Z}^+$).
However, there is no uniqueness of the decomposition of $f^n$ with the DEO family ${\Psi}_{k}^{+}$ ($k\in \mathbb{Z}$). A simple counter example can be found using the DEO family:
\begin{equation}\label{definitionOfeta}
 \eta_k(f) = 3(\partial_t f \partial_t^{k-1} f)-f\partial_t^k f , \qquad \forall k \in \mathbb{Z}
\end{equation}
Note that the derivative chain rule property is applied to this operator. One can verify:
 \begin{eqnarray}\label{Coeff2nnn}
\partial_t f^2 &=& 2 f \partial_t f \nonumber \\
\partial_t f^2 &=& \eta_1(f) \nonumber \\ 
\eta_1(f)&=&{\Psi}_{1}^{+}(f) \nonumber \\
& & \nonumber \\
\partial_t^2 f^2 &=& 2 (\partial_t f)^2 + 2 f \partial_t^2 f \nonumber \\
\partial_t^2 f^2 &=&  \partial_t {\eta}_{1}^{+}(f) \nonumber \\
\partial_t {\Psi}_{1}^{+}(f) &=& \partial_t {\eta}_{1}^{+}(f) \nonumber \\
\partial_t^2 f^2 &=& {\eta}_{2}^{+}(f) + {\eta}_{0}^{+}(\partial_t f) \nonumber \\
& & \nonumber \\
\partial_t^3 f^2 &=&  \partial_t^2 {\eta}_{1}^{+}(f) \nonumber \\
\partial_t^2 {\Psi}_{1}^{+}(f) &=& \partial_t^2 {\eta}_{1}^{+}(f) \nonumber \\
\partial_t^3 f^2 &=& {\eta}_{3}^{+}(f) +2 {\eta}_{1}^{+}(\partial_t f) + {\eta}_{-1}^{+}(\partial_t^2 f) \nonumber \\
& & \nonumber \\
\partial_t^4 f^2 &=&  \partial_t^3 {\eta}_{1}^{+}(f) \nonumber \\
\partial_t^3 {\Psi}_{1}^{+}(f) &=& \partial_t^3 {\eta}_{1}^{+}(f) \nonumber \\
\partial_t^4 f^2 &=& {\eta}_{4}^{+}(f) +3 {\eta}_{2}^{+}(\partial_t f) +3 {\eta}_{0}^{+}(\partial_t^2 f) + {\eta}_{-2}^{+}(\partial_t^3 f) 
\end{eqnarray}
\end{proof}
It is important to underline the fact that $(a_k^+)_{k\in\mathbb{Z}^+-\{0\}}$, $(A_k^+)_{k\in\mathbb{Z}^+-\{0\}}$ and $(B_k^+)_{k\in\mathbb{Z}^+-\{0\}}$ are not included in the definitions $1$ and $2$ as they do not follow the derivative chain rule.  In addition, $Definition$ $1$ can also be simplified fixing $l=2$ based on the definition of the energy operator in Equation \eqref{Psik+defdefef}.
\vspace{1.5em}
\newline $\bold{Theorem}$: for $f$ in $\mathbf{s}^{-}(\mathbb{R})$, the family of DEOs ${\Psi}_{k}^{-}$  and ${\Psi}_{k}^{+}$ ($k\in \mathbb{Z}$) decomposes uniquely the successive derivatives of the $n$-th power of $f$ for $n\in\mathbb{Z}^+$ and $n>1$. 
\begin{proof}
The proof is separated into three parts. The preliminary part justifies by induction the decomposition of the successive derivatives of the $n$-th power of $f$ for $n\in\mathbb{Z}^+$ and $n>1$ with ${\Psi}_{k}^{-}$  and ${\Psi}_{k}^{+}$ ($k\in \mathbb{Z}$). The proof is similar to the one in the Lemma, hence some parts are shorten to avoid long repetitions.  
\\The second section focuses on the existence of the decomposition as defined in Definition $1$ and Definition $2$. This is also shown by induction on the successive derivatives of $f^n$ ($n\in\mathbb{Z}^+$ and $n>1$) reusing the examples of different families of energy operator previously seen. Finally, the uniqueness is detailed in the last part.
\vspace{1.5em}
\newline A- $\bold{Preliminary}$
\vspace{1.5em}
\newline Recall the formula of ${\Psi}_{1}^{-}(f)$:
\begin{equation}
{\Psi}_{1}^{-}(f) = \dot{f}f - f \dot{f} =0
\end{equation}
One can see from the definition of ${\Psi}_{1}^{-}$, that this DEO projects any smooth real-valued function $f$ (in $\mathbf{s}^{-}(\mathbb{R})$) onto the null space (${\Psi}_{1}^{-}(f)$: $\mathbb{R}\times\mathbb{R}$ $\rightarrow$ $0$).
For any $f$ in $\mathbf{s}^{-}(\mathbb{R})$,
\vspace{1.5em}
\newline $\bold{Case}$ $n=2$:
\vspace{1.5em}
\begin{eqnarray}\label{equationPsikm01}
\partial_t f^2 &=& f \partial_t f + f\partial_t f + f\partial_t f - f\partial_t f \nonumber \\
\partial_t f^2 &=& {\Psi}_{1}^{+}(f) + {\Psi}_{1}^{-}(f)\nonumber \\
& & \nonumber \\
\partial_t^2 f^2 &=& 2 (\partial_t f)^2 + 2 f \partial_t^2 f \nonumber \\
\partial_t^2 f^2 &=& \partial_t ({\Psi}_{1}^{+}(f) +{\Psi}_{1}^{-}(f)) \nonumber \\
\partial_t^2 f^2 &=& {\Psi}_{2}^{+}(f) + {\Psi}_{0}^{+}(\partial_t f) + \nonumber \\
                 &&   {\Psi}_{2}^{-}(f) + {\Psi}_{0}^{-}(\partial_t f) \nonumber \\
& & \nonumber \\
\partial_t^3 f^2 &=& ({\Psi}_{3}^{+}(f) + {\Psi}_{3}^{-}(f)) +2 ({\Psi}_{1}^{+}(\partial_t f)+{\Psi}_{1}^{-}(\partial_t f)) + \nonumber \\
& & ({\Psi}_{-1}^{+}(\partial_t^2 f) +{\Psi}_{-1}^{-}(\partial_t^2 f)) \nonumber \\
\end{eqnarray}
In this development, there is a symmetry with the proof of the previous lemma (e.g., Equation \eqref{Coeff2n}). Thus, this allows to shorten the proof by induction of the decomposition of the successive derivatives of the $n$-th power of $f$. Therefore, a necessary step as shown in the demonstration of the existence in the Lemma, is the demonstration of a similar formula for $\partial_t^k ({\Psi}_{1}^{-}(f))$ (e.g., Equation \eqref{coefficientsAPgeneralbis2}). It ends with the generalization for $n=p$.
%
%
\\ From Equations \eqref{equationPsikm01}, one can define $a_s^-(f)$ in the same way that $a_s^+(f)$ was defined in Equations \eqref{Aocorf1N2} as:
\begin{equation}\label{Aocorf1N3}
a_s^-(f) = \partial_t^{s-1} \Psi_{1}^{-} (f), \qquad \forall s\in \mathbb{Z}^+-\{0\}
\end{equation}
%
%
%
With the property of the derivative chain rule in Equation \eqref{derivative chain rules0111} and Equation \eqref{Aocorf1N2}, it is easy to calculate the first terms of the DEO family $a_s^-(f)$ such as :
\begin{eqnarray}\label{Coeff2nmoins}
%
a_1^-(f)={\Psi}_{1}^{-}(f) &=& 0 \nonumber \\
& & \nonumber \\
a_2^-(f)=\partial_t {\Psi}_{1}^{-}(f) &=& {\Psi}_{2}^{-}(f) + {\Psi}_{0}^{-}(\partial_t f) \nonumber \\
{\Psi}_{2}^{-}(f) &=& -{\Psi}_{0}^{-}(\partial_t f) \nonumber \\
& & \nonumber \\
a_3^-(f)=\partial_t^2 {\Psi}_{1}^{-}(f) &=& {\Psi}_{3}^{-}(f) +2 {\Psi}_{1}^{-}(\partial_t f) + {\Psi}_{-1}^{-}(\partial_t^2 f) \nonumber \\
{\Psi}_{3}^{-}(f) &=& -{\Psi}_{-1}^{-}(\partial_t^2 f) \nonumber \\
& & \nonumber \\
a_4^-(f)=\partial_t^3 {\Psi}_{1}^{-}(f) &=& {\Psi}_{4}^{-}(f) +3 {\Psi}_{2}^{-}(\partial_t f) +3 {\Psi}_{0}^{-}(\partial_t^2 f) + {\Psi}_{-2}^{-}(\partial_t^3 f) \nonumber \\
{\Psi}_{4}^{-}(f) &=& -{\Psi}_{-2}^{-}(\partial_t^3 f)
\end{eqnarray}
The above equations are similar to the development in Equation \eqref{Coeff2n}. 
In addition, the family of DEO ${\Psi}_{k}^{-}$ ($k\in\mathbb{Z}$) has  the same derivative properties as ${\Psi}_{k}^{+}$. A similar equation can then be established for $a_s^-(f)$ following the development written in Equation \eqref{Coeff2n} as:
\begin{eqnarray}\label{coefficientsAPgeneral2bis}
a_s^-(f) &=& \sum_{k=0}^{s-1} \big(_{k}^{s-1} \big) \Psi_{2(k+1)-s}^{-}(\partial_t^{s-k-1}f), \forall s\in \mathbb{Z}^{+} -\{0\}
\end{eqnarray}
and, 
\begin{equation}
 a_s^-(f) = \partial_t^s \Psi_1^{-}(f) =0, \qquad \forall s \in \mathbb{Z}^+-\{0\}, \qquad \forall f\in\mathbf{s}^{-}(\mathbb{R})  
\end{equation}
This formula has just been checked for $s=\{1,2,3,4\}$ with the Equation \eqref{Coeff2nmoins}.
The generalization of the formula for $s=m$ is very similar to that described in the Equations \eqref{coefficientsAPgeneralbis} literally by changing $+$ and $-$ in the definition of the energy operator. It follows that the decomposition of the successive derivatives of $f^2$ is generalized for any $p$ in $\mathbb{Z}^+-\{0\}$ as:
\begin{eqnarray}\label{Coeff2nmoinspp}
\partial_t^m f^2 &=& \partial_t^{m-1} ({\Psi}_{1}^{+}(f) +{\Psi}_{1}^{-}(f)) \nonumber \\
                 &=& \sum_{k=0}^{m-1} \big(_{k}^{m-1} \big) \Psi_{2(k+1)-m}^{+}(\partial_t^{m-k-1}f) + \nonumber \\
                 & & \sum_{k=0}^{m-1} \big(_{k}^{m-1} \big) \Psi_{2(k+1)-m}^{-}(\partial_t^{m-k-1}f) \nonumber \\
\end{eqnarray}
In this case, $({\Psi}_{k}^{+})_{k\in\mathbb{Z}}$  and $({\Psi}_{k}^{-})_{k\in\mathbb{Z}}$ directly decompose $\partial_t^m f^2$. Note that the family $(a_k^-)_{k\in\mathbb{Z}^+-\{0\}}$ does not follow the derivative chain rule, and thus cannot be defined using definitions $1$ and $2$. In addition, the Equation \eqref{coefficientsAPgeneral2bis} can be easily extended for all $f$ in $\mathbf{S}^{-}(\mathbb{R})$ if we follow the induction proof in the Lemma (e.g., Equations \eqref{coefficientsAPgeneralbis}) as there is no restriction on the Images of the DEOs $({\Psi}_{k}^{-})_{k\in\mathbb{Z}}$.
\vspace{1.5em}
\\$\bold{Case}$ $n=p$ with $p>1$: 
\vspace{1.5em}
\\ Following the same step as in the proof of the previous lemma, let us define the families of DEO ${\theta}_{k}^{-}$ proportional to the family of DEO  ${\Psi}_{k}^{-}$  and ${\theta}_{k}^{+}$ proportional to the family of DEO  ${\Psi}_{k}^{+}$ ($k\in\mathbb{Z}$) with the assumption that they decompose the successive derivatives of ${f}^{p-1}$ as:
\begin{eqnarray}\label{psi+0008}
{\theta}_{k}^{+}(f) &=& \frac{(p-1)}{2} (\dot{f}{f}^{(k-1)} + f {f}^{(k)}) \nonumber \\
&= & \frac{(p-1)}{2} \Psi_k^+(f) \nonumber \\
{\theta}_{k}^{-}(f) &= & \frac{(p-1)}{2} (\dot{f}{f}^{(k-1)} - f {f}^{(k)}) \nonumber \\
&= &\frac{(p-1)}{2} \Psi_k^-(f)
\end{eqnarray}
%
Following the same development as in Equation \eqref{NONOVO}, one can define  $B_i^+(f) = \partial_t^{i-1} {\theta}_1^+(f)$ and $B_i^-(f) = \partial_t^{i-1} {\theta}_1^-(f)$ (with $i$ in $\mathbb{Z}^{+}-\{0\}$). Let us then write the first derivatives of ${f}^{p}$ such as:
\begin{eqnarray}\label{NONOVO2}
 \partial_t f^p &=& p f^{p-1} \partial_t f  \nonumber \\
  \partial_t f^p &=& \frac{p}{2} f^{p-2} (\Psi_k^+(f) +\Psi_k^-(f))  \nonumber \\
 \partial_t f^p &=& \frac{p}{p-1} ({\theta}_{1}^{+}(f)+{\theta}_{1}^{-}(f)) f^{p-2} \nonumber \\
\partial_t f^p &=& \frac{p}{p-1} (B_1^+(f) + B_1^-(f))  f^{p-2} \nonumber \\
& & \nonumber \\
%
\partial_t^2 f^p &=&  \frac{p}{p-1} (B_1^+(f)+B_1^-(f)) \partial_t f^{p-2}+ \frac{p}{p-1} (B_2^+(f)+B_2^-(f))  f^{p-2} )\nonumber \\
\end{eqnarray}
There is again a symmetry with the proof in the lemma. Following Equation \eqref{dtppfwithAps}, we can define the $s+1$-th derivative of $f^p$ using $B_{k+1}^-(f)$ and $B_{k+1}^+(f)$: 
\begin{equation}\label{dtppfwithApsd2}
\partial_t^{s+1} f^p = \sum_{k=0}^s \big(_{k}^s \big) \frac{p}{p-1}(B_{k+1}^-(f) +B_{k+1}^+(f)) \partial_t^{s-k} f^{p-2}, \qquad \forall s \in \mathbb{Z}^+
\end{equation}
This equation has just been checked for $s=\{0,1\}$. As the induction proof follows exactly the same development as in Equation \eqref{dtp3fwithAp02} by only adding $B_{k+1}^-(f)$ (with the same properties as $B_{k+1}^+(f)$ in Equation \eqref{NONOVO2}), it allows then to assume the generalization to the case $s+2$.
\newline Thus,  $(B_{k}^+)_{k\in\mathbb{Z}}$ and  $(B_{k}^-)_{k\in\mathbb{Z}}$ decompose the $s$-th derivative of $f^p$. From their definition, one can conclude that $(\theta_{k}^+)_{k\in\mathbb{Z}}$ and $(\theta_{k}^-)_{k\in\mathbb{Z}}$ decompose $\partial_t^sf^p$.
%
%
$({\theta}_{k}^{+})_{k\in\mathbb{Z}}$ and $({\theta}_{k}^{-})_{k\in\mathbb{Z}}$ are proportional to  $({\Psi}_{k}^{+})_{k\in\mathbb{Z}}$  and $({\Psi}_{k}^{-})_{k\in\mathbb{Z}}$ respectively. Finally, $({\Psi}_{k}^{+})_{k\in\mathbb{Z}}$  and $({\Psi}_{k}^{-})_{k\in\mathbb{Z}}$  decompose the successive derivatives of ${f}^{p}$ ($p \in \mathbb{Z}^+$, and $p> 1$).
\vspace{1.5em}
\newline B- $\bold{Existence}$ $\bold{of}$ $\bold{the}$ $\bold{Decomposition}$ 
\vspace{1.5em}
\newline The proof is also structured as an induction on the $n$-th power of $f$. The different cases revisit some families of operator defined in the previous sections of this work (e.g., proof of the Lemma).
%
\vspace{1.5em}
\\$\bold{Case}$ $n=2$: 
%
%
It was shown that the family of operator $(\eta_k)_{k\in\mathbb{Z}}$ (see definition in the proof of the Lemma), decomposes $\partial^s_t f^{2}$ ($s\in\mathbb{Z}^+-\{0\}$). As defined in Equation \eqref{definitionOfeta}, one can rewrite it as a sum of the DEO family $({\Psi}_{k}^{-})_{k\in\mathbb{Z}}$ and $({\Psi}_{k}^{+})_{k\in\mathbb{Z}}$ as:
\begin{equation}
\eta_k(f) = \Psi_k^+(f) + 2 \Psi_k^-(f), \qquad k \in \mathbb{Z}
\end {equation}
\vspace{1.5em}
\\$\bold{Case}$ $n=3$: 
%
%
In the lemma, the family of operator $(\Gamma_{k}^{+})_{k\in\mathbb{Z}}$ was defined in Equation \eqref{GammaOperatorfamily} and decomposes $\partial^s_t f^{3}$ ($s\in\mathbb{Z}^+-\{0\}$). One can rewrite the definition as:
\begin{eqnarray}\label{GammaOperatorfamilyPi}
%
\Gamma_{k}^{+} (f) &=& \frac{3}{2} (\dot{f}{f}^{(k-1)} + f {f}^{(k)}) \nonumber \\
\Gamma_{k}^{+} (f) &=& \frac{3}{2} \Psi_k^+(f) + 0 \Psi_k^-(f)
\end{eqnarray}
\vspace{1.5em}
\\$\bold{Case}$ $n=p$ with $p>1$: 
%
%
Previously, Equations \eqref{psi+0008} defined the operators ${\theta}_{k}^{+}$ and ${\theta}_{k}^{-}$ decomposing $\partial^{s}_t f^{p}$ ($s\in\mathbb{Z}^+-\{0\}$) and proportional to $({\Psi}_{k}^{-})_{k\in\mathbb{Z}}$ and $({\Psi}_{k}^{+})_{k\in\mathbb{Z}}$ as:
%
\begin{eqnarray}\label{psi+000100}
{\theta}_{k}^{+}(f) &=& \frac{(p-1)}{2} (\dot{f}{f}^{(k-1)} + f {f}^{(k)}) \nonumber \\
&= & \frac{(p-1)}{2} \Psi_k^+(f) \nonumber \\
{\theta}_{k}^{-}(f) &= & \frac{(p-1)}{2} (\dot{f}{f}^{(k-1)} - f {f}^{(k)}) \nonumber \\
&= &\frac{(p-1)}{2} \Psi_k^-(f)
\end{eqnarray}
By induction it was also shown in the same section that ${\theta}_{k}^{-}$ and ${\theta}_{k}^{+}$ decompose $\partial^{s}_t f^{p+1}$ ($s\in\mathbb{Z}^+-\{0\}$). It is then possible to conclude the existence of the decomposition of any operator by using $({\Psi}_{k}^{-})_{k\in\mathbb{Z}}$ and $({\Psi}_{k}^{+})_{k\in\mathbb{Z}}$.
%
%
%
%
\vspace{1.5em}
\newline C-  $\bold{About}$ $\bold{Uniqueness}$ $\bold{of}$ $\bold{the}$  $\bold{Decomposition}$ 
\vspace{1.5em}
%
\\ With the previous section, let us show by induction the uniqueness of the decomposition of any family of operators decomposing $\partial_t^s f^n$ ($s$ in $\mathbb{Z}^+$, $n$ in $\mathbb{Z}^+$ and $n > 1$). The induction is focused on the proof of the uniqueness of the decomposition of a family of operator $({S}_{k})_{k\in\mathbb{Z}}$ (following the derivative chain rule property) by $({\Psi}_{k}^+)_{k\in\mathbb{Z}}$ and $({\Psi}_{k}^-)_{k\in\mathbb{Z}}$. In other words, the induction is on the $k$-th order of the operator. 
%
%
\vspace{1.5em}
\\$\bold{Case}$ $k=2$: 
For $f$ in $\mathbf{s}^{-}(\mathbb{R})$ and $n$ in $\mathbb{Z}^+$ and $n > 1$, one can assume that $(\alpha_1,\alpha_2,\beta_1,\beta_2)$ exist in $\mathbb{R}^4$ such as:
%
%
%
%
\begin{eqnarray}\label{eqrefpartial0102}
\partial_t^s f^n &=& \partial_t^{s-1} (f^{n-2} S_1(f))  \nonumber \\
\partial_t^s f^n &=&  \partial_t^{s-1} ( f^{n-2} (\alpha_1 \Psi^{+}_1(f) +\alpha_2 \Psi^{-}_1(f))) \nonumber \\
\partial_t^s f^n &=& \partial_t^{s-1} ( f^{n-2} (\beta_1 \Psi^{+}_1(f) +\beta_2 \Psi^{-}_1(f))) \nonumber \\
\end{eqnarray} 
%
%
\\ As  with the operator family $({S}_{k})_{k\in\mathbb{Z}}$ follows the derivative chain rule property:
%
\begin{eqnarray}\label{equationalphabeta1}
\partial_t S_1(f) &=& S_2(f) +S_0(\partial_t f)\nonumber \\
\partial_t S_1(f) &=&\alpha_1 \partial_t \Psi^{+}_1(f) +\alpha_2 \partial_t \Psi^{-}_1(f) \nonumber \\
\partial_t S_1(f) &=&\alpha_1 (\Psi^{+}_2(f) +\Psi^{+}_0(\partial_t f)) +\alpha_2 ( \Psi^{-}_2(f) +\Psi^{-}_0(\partial_t f)) \nonumber \\
\end{eqnarray} 
And then,
\begin{eqnarray}\label{equationalphabeta1b}
S_2(f) &=& \alpha_1 \Psi^{+}_2(f) +\alpha_2 \Psi^{-}_2(f) \nonumber \\
S_2(f) &=& \beta_1 \Psi^{+}_2(f) +\beta_2 \Psi^{-}_2(f) \nonumber \\
(\alpha_1-\beta_1) \Psi^{+}_2(f) +(\alpha_2-\beta_2) \Psi^{-}_2(f) &=& 0
\end{eqnarray} 
As $Im(\Psi^{+}_2)$ and $Im(\Psi^{-}_2)$ are not reduced to $\{0\}$ by definition, it follows that $\alpha_1=\beta_1$ and $\alpha_2=\beta_2$. Note that it is not possible to do this simple check for $k=1$ as $Im(\Psi^{-}_1)=\{0\}$.
\vspace{1.5em}
\\$\bold{Case}$ $k=p$: 
Let us assume the uniqueness of the decomposition for $k=p-1$ (with $k\neq1$). For $k=p$, following Equation \eqref{equationalphabeta1}:
\begin{eqnarray}\label{equationalphabeta1kn}
\partial_t S_{p-1}(f) &=& S_p(f) +S_{p-2}(\partial_t f)\nonumber \\
\partial_t S_{p-1}(f) &=&\alpha_1 \partial_t \Psi^{+}_{p-1}(f) +\alpha_2 \partial_t \Psi^{-}_{p-1}(f) \nonumber \\
\partial_t S_{p-1}(f) &=&\alpha_1 (\Psi^{+}_p(f) +\Psi^{+}_{p-2}(\partial_t f)) +\alpha_2 ( \Psi^{-}_p(f) +\Psi^{+}_{p-2}(\partial_t f)) \nonumber \\
\end{eqnarray} 
And then,
\begin{eqnarray}\label{equationalphabeta1knb}
S_p(f) &=& \alpha_1 \Psi^{+}_p(f) +\alpha_2 \Psi^{-}_p(f) \nonumber \\
S_p(f) &=& \beta_1 \Psi^{+}_p(f) +\beta_2 \Psi^{-}_p(f) \nonumber \\
(\alpha_1-\beta_1) \Psi^{+}_p(f) +(\alpha_2-\beta_2) \Psi^{-}_p(f) &=& 0
\end{eqnarray} 
By definition for $p\neq1$, $Im(\Psi^{+}_p)$ and $Im(\Psi^{-}_p)$ are not reduced to $\{0\}$, and it follows that $\alpha_1=\beta_1$ and $\alpha_2=\beta_2$. 
\vspace{1.5em}
\\$\bold{Special}$ $\bold{Case}$ $k=1$: 
To complete the proof with the assumption that  $\alpha_1=\beta_1$ and $\alpha_2=\beta_2$ for $k \in\mathbb{Z}$ and $k \neq 1$ , the special case $k=1$ can be solved as:
\begin{eqnarray}\label{equationalphabeta1c}
\partial_t(\alpha_1 \Psi^{+}_1(f)) &=&\alpha_1(\Psi^{+}_2(f) + \Psi^{+}_0(\partial_t f)) \nonumber \\
&=&\beta_1(\Psi^{+}_2(f) + \Psi^{+}_0(\partial_t f)) \nonumber \\
&=&\partial_t(\beta_1\Psi^{+}_1(f))\nonumber \\
& & \nonumber \\
\partial_t(\alpha_2 \Psi^{-}_1(f)) &=&\alpha_2(\Psi^{-}_2(f) + \Psi^{-}_0(\partial_t f)) \nonumber \\
&=&\beta_2(\Psi^{-}_2(f) + \Psi^{-}_0(\partial_t f)) \nonumber \\
&=&\partial_t(\beta_2\Psi^{-}_1(f))\nonumber \\
\nonumber\\
\end{eqnarray} 
Thus, $\alpha_1=\beta_1$ and $\alpha_2=\beta_2$. Let us finish the proof with a remark for the case $\partial_t f^2$. Rewriting Equation \eqref{eqrefpartial0102} as:
\begin{eqnarray}
\partial_t f^2 &=& S_1(f)  \nonumber \\
\partial_t f^2 &=& \alpha_1 \Psi^{+}_1(f) +\alpha_2 \Psi^{-}_1(f) \nonumber \\               
\partial_t f^2 &=& \beta_1 \Psi^{+}_1(f) +\beta_2  \Psi^{-}_1(f) \nonumber \\
&=&  \Psi^{+}_1(f) 
\end{eqnarray} 
Thus, we conclude that in this case, $\alpha_1 = \beta_1 =1$.
%
\end{proof}
Note that in \cite{Boudraa et al.2009}, the authors based their work on the energy operator defined as $LP_2(f,g)=\partial_t f \partial_t g - f \partial_t^2 g$ ($LP_2$ : $\mathbb{R}^2$ $\rightarrow$ $\mathbb{R}$) which was then generalized to the complex set ($LP_2^{\mathbb{C}}(f,g)$). In the development of their work, they found a similar type of formula as found in Equation \eqref{coefficientsAPgeneral}, but restricted to the definition of their energy operator. In addition, one can underline that the generalization of the decomposition of the successive derivatives of the $n$-th power of $f$ with the DEOs ${(\Psi^{+}_k)}_{k\in \mathbb{Z}}$ and ${(\Psi^{-}_k)}_{k\in \mathbb{Z}}$ (e.g., Equations \eqref{coefficientsAPgeneral} and \eqref{Coeff2nmoinspp}) follows the general Leibniz derivative rules \cite{BruceWest}. 
%
%
%
%
\vspace{0.5em}
\newline $\bold{Discussion}$ $n<-1$: This case focuses on the decomposition using the DEO family of the quotient of the function $f$:($\mathbb{R}\rightarrow\mathbb{R}$) defined as:
\begin{equation}
\forall f \in \mathbf{S}^{-}(\mathbb{R}), \qquad \forall t \in \mathbb{R}, \qquad f(t)\neq 0, \qquad \forall n \in \mathbb{Z}^+, n>1, \frac{1}{f^n}
\end{equation}
They are just a particular case of $f^n$ ($n>1$, $n\in \mathbb{Z}^+$). Using an intermediary function, $h$ such as $h = \frac{1}{f}$, the problem of decomposing $\partial_t^s f^{-n}$ ($s \in \mathbb{Z}^+$) is equivalent to resolving $\partial_t^s h^{n}$, which has been demonstrated in the Lemma and Theorem.
\vspace{1.5em}
\newline $\bold{Discussion}$ $n=1$: This case does not make sense to decompose $f$ and its $k$-th derivative ($k\in\mathbb{Z}^+$) as a sum of energy operators based on the general definition of the DEO from \cite{Maragos1995}: the $k$-th order DEO is the cross energy between a function and its $k-1$ derivatives. However, one can use a general formula in the special case:
\begin{equation}
\forall t\in \mathbb{R} , \qquad f\in \mathbf{S}^{-}(\mathbb{R}),\qquad f(t) \neq 0 \nonumber
\end{equation}
\begin{eqnarray}\label{discussion2a}
\partial_t^k f &=& \partial_t^k \big(\frac{f^3}{f^2}\big) \nonumber \\
k=1, \qquad \partial_t f &=& f^{-2} \partial_t f^3 + f^3 \partial_t f^{-2} \nonumber \\
k=2, \qquad \partial_t^2 f &=& 2 \partial_t f^{-2} \partial_t f^3 + f^3 \partial_t^2 f^{-2} + f^{-2} \partial_t^2 f^3
\end{eqnarray} 
The example for $k =\{1,2\}$ in Equation \eqref{discussion2a} shows that $\partial_t f$  can be decomposed into a product of successive derivatives of $f^{3}$ and $f^{-2}$. Those derivatives can be decomposed into a sum of DEOs based on the development of the proof of the Lemma and the previous discussion (for $n<-1$ ).
\vspace{0.5em}
%
%
%
\newline Note that if we restrict $f$ in $\mathbf{S}^{-}(\mathbb{R})$ and to be expandable in Taylor-Series in an interval $[a,b]$ in $\mathbb{R}$ such as:
\begin{equation}
f(t) = \Sigma_{k=0}^{\infty} (\partial_t^k f(t_0)) \frac{(t-t_0)^k}{k!}, \qquad \forall t \in [a,b], \qquad t_0 \in [a,b] 
\end{equation}
One can then decompose the coefficients of the Taylor-Series of $f$ as a sum of DEO via the method shown in the Lemma and Theorem.
%
%
%
%
%
%
%
%
%
\section{Some Properties of the Images and Kernels of the DEO Families}\label{SectionProperties}
In this part, we study the relationships between $Im(\Psi^{+}_k)$ and $Im(\Psi^{-}_k)$, and $Ker(\Psi^{+}_k)$ and $Ker(\Psi^{-}_k)$  via the demonstration of the following properties.
\\$\mathbf{Properties}$ $1$: for $k$ in $\mathbb{Z}$ and $f$ in $\mathbf{S}^{-}(\mathbb{R})$,
\begin{eqnarray}\label{Property01}
a)  \Psi^{+}_k(f) &=& \Psi^{+}_{-k+2}(\partial_t^{k-1}f) \nonumber \\
b) \Psi^{-}_k(f) &=& -\Psi^{-}_{-k+2}(\partial_t^{k-1}f) \nonumber \\
c) \Psi^{+}_k(f)+\Psi^{-}_k(f) &=& \Psi^{+}_{-k+2}(\partial_t^{k-1}f) -\Psi^{-}_{-k+2}(\partial_t^{k-1}f) \nonumber \\ 
d) Im(\Psi^{+}_k) \bigcap Im(\Psi^{-}_k) & \subseteq & Im(\Psi^{+}_k - \Psi^{-}_k) \nonumber \\
e) Im(\Psi^{+}_k) \bigcap Im(-\Psi^{-}_k) & \subseteq & Im(\Psi^{+}_k + \Psi^{-}_k) \nonumber \\
\end{eqnarray}
\begin{proof}
To show Equation (\ref{Property01}-a), the definition of $\Psi^{+}_{-k+2}(\partial_t^{k-1}f)$ is:
\begin{equation}
\Psi^{+}_{-k+2}(f) = \partial_t f \partial_t^{-k+1}f + f \partial_t^{-k+2}f
\end{equation}
and then,
\begin{eqnarray}\label{equationproperty01b}
\Psi^{+}_{-k+2}(\partial_t^{k-1}f) &=& \partial_t^k f f + \partial_t^{k-1}f \partial_t f \nonumber \\
\Psi^{+}_{-k+2}(\partial_t^{k-1}f) &=& \Psi^{+}_k(f) 
\end{eqnarray}
This last equation then shows the first assertion $\Psi^{+}_k(f) = \Psi^{+}_{-k+2}(\partial_t^{k-1}f)$. Following the same development as in Equation \eqref{equationproperty01b}, one can derive directly Equation (\ref{Property01}-b) as:
\begin{eqnarray}\label{equationproperty01c}
\Psi^{-}_{-k+2}(\partial_t^{k-1}f) &=& \partial_t^k f f - \partial_t^{k-1}f \partial_t f \nonumber \\
\Psi^{-}_{-k+2}(\partial_t^{k-1}f) &=& -\Psi^{-}_k(f) 
\end{eqnarray}
Equation (\ref{Property01}-c) is a consequence of  Equation (\ref{Property01}-a) and Equation (\ref{Property01}-b) by simply replacing  $\Psi^{+}_k(f)$ and $\Psi^{-}_k(f)$ by respectively $\Psi^{+}_{-k+2}(\partial_t^{k-1}f)$ and $-\Psi^{-}_{-k+2}(\partial_t^{k-1}f)$. 
\newline Let us write the definition of $Im(\Psi^{+}_k) \bigcap Im(\Psi^{-}_k) $ and $Im(\Psi^{+}_k - \Psi^{-}_k)$ such as:
\begin{eqnarray}
 Im(\Psi^{+}_k) \bigcap Im(\Psi^{-}_k) &=& \{\Psi^{+}_k(f) = \Psi^{-}_k(f) |f \in \mathbf{S}^{-}(\mathbb{R})\} \nonumber \\ 
                                                             &=& \{ 2 f\partial_t^k f =0 |f \in \mathbf{S}^{-}(\mathbb{R})\} \nonumber \\
Im(\Psi^{+}_k - \Psi^{-}_k) &=& \{\Psi^{+}_k(f) - \Psi^{-}_k(f) |f \in \mathbf{S}^{-}(\mathbb{R})\} \nonumber \\
                                           &=& \{ 2 f\partial_t^k f |f \in \mathbf{S}^{-}(\mathbb{R})\} 
\end{eqnarray}
$Im(\Psi^{+}_k - \Psi^{-}_k)$ is non-empty as it contains $0$ when $f$ is the null function. Thus, the above definitions show the inclusion of the subset $Im(\Psi^{+}_k) \bigcap Im(\Psi^{-}_k)$ into $Im(\Psi^{+}_k - \Psi^{-}_k)$. In the same way, one can show Equation (\ref{Property01}-e) through the definition of each subset.
\begin{eqnarray}
 Im(\Psi^{+}_k) \bigcap Im(-\Psi^{-}_k) &=& \{\Psi^{+}_k(f) +\Psi^{-}_k(f) =0 |f \in \mathbf{S}^{-}(\mathbb{R})\} \nonumber \\ 
Im(\Psi^{+}_k + \Psi^{-}_k) &=& \{\Psi^{+}_k(f) + \Psi^{-}_k(f) |f \in \mathbf{S}^{-}(\mathbb{R})\} \nonumber \\                              
\end{eqnarray}
As before, $0$ is included in $Im(\Psi^{+}_k + \Psi^{-}_k)$ (e.g. $\Psi^{+}_k (0) = \Psi^{-}_k (0) =0$, for $k$ in $\mathbb{Z}$). Thus, $Im(\Psi^{+}_k) \bigcap Im(-\Psi^{-}_k) \subset Im(\Psi^{+}_k + \Psi^{-}_k)$.
Note that Equations (\ref{Property01}-a), (\ref{Property01}-b) and (\ref{Property01}-c) are important as they are directly linked to $Im(\partial_t^p\Psi^{+}_1)$ and $Im(\partial_t^p\Psi^{-}_1)$ ($p$ in $\mathbb{Z}-\{0\}$).
%
%
\end{proof}
Similar properties can be shown from the definitions of $Ker(\Psi^{+}_k)$ and $Ker(\Psi^{-}_k)$:
\\$\mathbf{Properties}$ $2$: for $k$ in $\mathbb{Z}$ and $f$ in $\mathbf{S}^{-}(\mathbb{R})$,
\begin{eqnarray}\label{Property02Kernel}
a) Ker(\Psi^{+}_k) \bigcap Ker(\Psi^{-}_k) & \subseteq & Ker(\Psi^{+}_k - \Psi^{-}_k) \nonumber \\
b) Ker(\Psi^{+}_k) \bigcap Ker(-\Psi^{-}_k) & \subseteq & Ker(\Psi^{+}_k + \Psi^{-}_k) \nonumber \\
c) Ker(-\Psi^{+}_k) &=& Ker(\Psi^{+}_k) \nonumber \\
d) Ker(-\Psi^{-}_k) &=& Ker(\Psi^{-}_k) \nonumber \\ 
\end{eqnarray}
%
\begin{proof}
%
%
%
%
%
%
%
The demonstration follows the definition of the kernels in the same way as in the previous properties:
\begin{eqnarray} \label{defKer01}
Ker(\Psi^{+}_k) \bigcap Ker(\Psi^{-}_k) &=& \{f\in \mathbf{S}^{-}(\mathbb{R})| \qquad \Psi^{-}_k(f)=\Psi^{+}_k(f) =0\} \nonumber \\
\end{eqnarray}
\begin{eqnarray}\label{defKer02}
Ker(\Psi^{+}_k - \Psi^{-}_k) &=& \{f\in \mathbf{S}^{-}(\mathbb{R})| \Psi^{+}_k(f) -\Psi^{-}_k(f) =0\} \nonumber \\
\end{eqnarray}
And, 
\begin{eqnarray} \label{defKer01b}
Ker(\Psi^{+}_k) \bigcap Ker(-\Psi^{-}_k)    &=& \{f\in \mathbf{S}^{-}(\mathbb{R})| \qquad \Psi^{+}_k(f)=-\Psi^{-}_k(f) =0\} \nonumber\\
\end{eqnarray}
\begin{eqnarray}\label{defKer02b}
Ker(\Psi^{+}_k + \Psi^{-}_k) &=& \{f\in \mathbf{S}^{-}(\mathbb{R})| \Psi^{+}_k(f) +\Psi^{-}_k(f) =0\}\nonumber \\
\end{eqnarray}
It is necessary to underline that $Ker(\Psi^{+}_k)$, $Ker(\Psi^{-}_k)$ , $Ker(\Psi^{+}_k + \Psi^{-}_k)$ and $Ker(\Psi^{+}_k - \Psi^{-}_k)$ are all non empty sets as the null function is included in all of them.  This then demonstrates Equations (\ref{Property02Kernel}-a) and (\ref{Property02Kernel}-b).
%
\newline Similarly we have,
\begin{eqnarray}\label{defKer03}
Ker(\Psi^{+}_k ) &=& \{f\in \mathbf{S}^{-}(\mathbb{R})| \Psi^{+}_k(f) = -\Psi^{+}_k(f) =0\} \nonumber \\
Ker(\Psi^{-}_k ) &=& \{f\in \mathbf{S}^{-}(\mathbb{R})| \Psi^{-}_k(f) = -\Psi^{-}_k(f) =0\} \nonumber \\
\end{eqnarray}
This directly shows Equations  (\ref{Property02Kernel}-c) and (\ref{Property02Kernel}-d). 
%
\end{proof}
Following the remarks at the end of Properties $1$, we can rewrite the formulas in Equations \eqref{coefficientsAPgeneral} based on the relationship of the images. With the Equation (\ref{Property01}-a), the formula \eqref{coefficientsAPgeneral} can be rewritten as:
\begin{eqnarray}\label{Coeff2nRewritten}
\partial_t {\Psi}_{1}^{+}(f) &=& {\Psi}_{2}^{+}(f) + {\Psi}_{0}^{+}(\partial_t f) \nonumber \\
\partial_t {\Psi}_{1}^{+}(f) &=& 2{\Psi}_{2}^{+}(f) \nonumber \\
\partial_t^2 {\Psi}_{1}^{+}(f) &=& {\Psi}_{3}^{+}(f) + 2{\Psi}_{1}^{+}(\partial_t f) + {\Psi}_{-1}^{+}(\partial_t^2 f) \nonumber \\
&=& 2{\Psi}_{3}^{+}(f) + 2{\Psi}_{1}^{+}(\partial_t f) \nonumber \\
\end{eqnarray}
%
%
%
Thus for $p$ in $\mathbb{Z}^{+}-\{0\}$, $k$ in $\mathbb{Z}$ and $s$ in in $\mathbb{Z}^+-\{0\}$,
\begin{equation*}\label{coefficientsAPgeneralRewrittenA}
a_p^+(f) =\left\{
 \begin{array}{rl}
 \Psi_{1}^{+}(f),&  p=1 \nonumber \\
2 \sum_{k > \frac{2s-1}{2}}^{2s-1} \big(_{k}^{2s-1} \big)\Psi_{2(k+1)-2s}^{+}(\partial_t^{2s-k-1}f), &  p=2s  \\
2 \sum_{k > s}^{2s} \big(_{k}^{2s} \big) \Psi_{2(k+1)-2s-1}^{+}(\partial_t^{2s-k}f) + \big(_{s}^{2s} \big) \Psi_{1}^{+}(\partial_t^{s}f),&  p=2s+1 \\
\end{array} \right .
\end{equation*}
Consequently, the successive derivatives of $f^n$ ($n$ in $\mathbb{Z}^{+}-\{0\}$) can be decomposed with the DEO family $(\Psi_{k}^{+})$ with $k$ in $\mathbb{Z}^+-\{0\}$. 
\newline In the same way following Equation (\ref{Property01}-b), Equation \eqref{coefficientsAPgeneral2bis} is simplified to:
\begin{equation}\label{coefficientsAPgeneralRewrittenA}
a_p^-(f) = 0 \qquad \forall  p \in \mathbb{Z}^{+}-\{0\}
\end{equation}
\section{Application to the Energy Function}
The energy function $\mathcal{E}$ considered in this section is the one previously defined in Equation \eqref{EnergyfunDefine02}.
$f_1$ and $f_1^2$ are considered to be in $\mathbf{S}^{-}(\mathbb{R})$, analytic and with a finite energy. 
Note that the choice of this example is also based on the discussion (case $n=1$) after the Lemma. 
%
%
\\ One can develop the energy function in Taylor-Series on the interval of definition $[0,\tau]$ ($\tau$ in $\mathbb{R}$) for a nominated $\tau_0$ in the defined interval,
\begin{eqnarray}\label{EquationEfunction}
\mathcal{E}(f_1(\tau) ) &=& \mathcal{E}(f_1(\tau_0)) + \sum_{k=1}^\infty \partial_t^k \mathcal{E}(f_1(\tau_0)) \frac{(\tau-\tau_0)^k}{k!} \nonumber \\
                                 &=& \mathcal{E}(f_1(\tau_0)) + \sum_{k=1}^\infty \partial_t^{k-1} f_1^2(\tau_0) \frac{(\tau-\tau_0)^k}{k!} \nonumber \\
                                  &=& \mathcal{E}(f_1(\tau_0)) + \sum_{k=1}^\infty \partial_t^{k-1} f_1^2(\tau_0) \frac{(\tau-\tau_0)^k}{k!} \nonumber \\
                                  &=& \mathcal{E}(f_1(\tau_0)) +f_1^2(\tau_0) (\tau-\tau_0) + \sum_{k=2}^\infty \partial_t^{k-2} (\Psi_1^+(f_1)(\tau_0)   \nonumber \\
& &+ \Psi_1^-(f_1)(\tau_0))\frac{(\tau-\tau_0)^k}{k!} \nonumber \\
&=& \mathcal{E}(f_1(\tau_0)) +f_1^2(\tau_0) (\tau-\tau_0) + \sum_{k=2}^\infty \partial_t^{k-2} \Psi_1^+(f_1)(\tau_0) \frac{(\tau-\tau_0)^k}{k!} \nonumber \\
\end{eqnarray}
In the special case that the series is absolutely convergent, one can write for $p$ in $\mathbb{Z}^+-\{0\}$:
\begin{eqnarray}
\big|  \frac{\partial_t^{p+1} \mathcal{E}(f_1(\tau_0)) }{\partial_t^{p} \mathcal{E}(f_1(\tau_0))} \big| \big| \frac{(\tau-\tau_0)}{p+1} \big|&<& 1 \nonumber \\
\end{eqnarray}
or for  $p$ in $\mathbb{Z}^+-\{0,1\}$
\begin{eqnarray}\label{equationenergyseries1}
\big| \frac{\partial_t^{p+1} (\Psi_1^+(f_1)(\tau_0) }{\partial_t^{p} \Psi_1^+(f_1)(\tau_0) } \big| &<& \big| \frac{p+1}{(\tau-\tau_0)}\big|
\end{eqnarray}
%
%
%
\newline Let us take an example with the periodic function :
\begin{equation}\label{funGt}
g(t) = A  \cos(t), \qquad t \in [-\pi,\pi], \qquad A \in \mathbb{R}
\end{equation}
We are interested in the development in Taylor-Series of the energy of $g$ following Equation \eqref{EquationEfunction} with $\tau$ and $\tau_0$ in $[-\pi,\pi]$  such that:
\begin{eqnarray}
\mathcal{E}(g(\tau))&=& \mathcal{E}(g(\tau_0)) +g^2(\tau_0) (\tau-\tau_0) + \sum_{k=2}^\infty \partial_t^{k-2} \Psi_1^+(g)(\tau_0) \frac{(\tau-\tau_0)^k}{k!} \nonumber \\
\end{eqnarray}
In addition, with the definition of $\Psi_1^+$ (e.g. Equation \eqref{Psik+defdefef}) one can write:
\begin{eqnarray}\label{cossinexample}
\Psi_1^+(g(t)) &=& -2A \cos(t) \sin(t)\nonumber \\
\partial_t \Psi_1^+(g(t)) &=& -2A (\cos^2(t) - \sin^2(t)) \nonumber \\
                          &=& 2A (2\sin^2(t)-1) \nonumber \\
\partial_t^2 \Psi_1^+(g(t)) &=& 8A \sin(t) \cos(t) \nonumber \\ 
\partial_t^3 \Psi_1^+(g(t)) &=& -8A (2\sin^2(t)-1) \nonumber \\ 
\end{eqnarray}
One can deduce a general formula for the derivatives of $g$ from the above equations as:
%
\begin{equation}
\forall k \in \mathbb{Z}^+,\left\{
\begin{array}{rcl}
 \partial_t^{2k+1} \Psi_1^+(g(t)) &=& (-1)^{k+1} 2^{2k+1} A (\cos^2(t) - \sin^2(t))\\
 \partial_t^{2k} \Psi_1^+(g(t)) &=& (-1)^{k+1} 2^{2k+1} A (\cos(t) \sin(t))\\
\end{array} \right.
\end{equation}
%
From those equations and the general trigonometric properties of the functions $cosines$ and $sines$, the upper bound of $\partial_t^{p} \Psi_1^+(g(t))$ is:
\begin{equation}
\forall k \in \mathbb{Z}^+,\left\{
\begin{array}{rcl}
 \partial_t^{2k+1} \Psi_1^+(g(t)) &\leq& | 2^{2k+1} A | \\
 \partial_t^{2k} \Psi_1^+(g(t)) &\leq& | 2^{2k+1} A | \\
\end{array} \right.
\label{example02eq}
\end{equation}
%
Let us now examine the convergence of the Taylor-Series of $\mathcal{E}(g)$ using the ratio test (see \cite{Kreyszig}) as described in Equation \eqref{equationenergyseries1}. One can write for $p$ in $\mathbb{Z}^+-\{0,1\}$:
\begin{eqnarray}
lim_{p\rightarrow + \infty}\big| \frac{\partial_t^{p+1} (\Psi_1^+(g)(\tau_0) }{\partial_t^{p} \Psi_1^+(g)(\tau_0) } \big| &=& \nonumber \\ 
lim_{p\rightarrow + \infty}\big| 2^2 \frac{(\tau-\tau_0)}{p+1}\big| &=& 0 
\end{eqnarray}
Thus, the Taylor-Series of $\mathcal{E}(g)$ is absolutely convergent. Moreover, we can find the distances between the successive derivatives of $\partial_t^{p} \Psi_1^+(g)(\tau_0)$ with $\tau$ and $\tau_0$ in $[-\pi,\pi]$.
\begin{eqnarray}
|\partial_t^{p+1} (\Psi_1^+(g)(\tau_0)| - |\partial_t^{p} (\Psi_1^+(g)(\tau_0)| &=& \big( 2^2 \frac{|\tau-\tau_0|}{p+1} -1 \big) |\partial_t^{p} (\Psi_1^+(g)(\tau_0)| \nonumber \\
|\partial_t^{p+1} (\Psi_1^+(g)(\tau_0)| - |\partial_t^{2} (\Psi_1^+(g)(\tau_0)| & =& \sum_{i=2}^{p+1} \big( 2^2 \frac{|\tau-\tau_0|}{i+1} -1 \big) |\partial_t^{i} (\Psi_1^+(g)(\tau_0)|\nonumber \\
\end{eqnarray}
With Equations \eqref{example02eq}, the distance can be upper bounded as:
\begin{equation}
|\partial_t^{p+1} (\Psi_1^+(g)(\tau_0)| - |\partial_t^{2} (\Psi_1^+(g)(\tau_0)| \left\{
\begin{array}{rl}
\leq & \sum_{k=1}^{(p+1)/2} \big( 2^2 \frac{|\tau-\tau_0|}{2 k+1} -1 \big) | 2^{2k+1} A | , \qquad p =2k\\
\leq & \sum_{k=1}^{p/2} \big( 2^2 \frac{|\tau-\tau_0|}{2 k+2} -1 \big) | 2^{2k+1} A | , \qquad p= 2k+1
\end{array} \right.
\end{equation}
Note that the distance for the case $p =\{0,1\}$ should be calculated with the equations in \eqref{cossinexample}.
\section{Conclusions}
This work defined two families of DEOs ${\Psi}_{k}^{+}$ and ${\Psi}_{k}^{-}$($k \in \mathbb{Z}$).
A lemma and theorem were developed to decompose any function $f$ in $\mathbf{S}^{-}(\mathbb{R})$ and its $n$-th power ($n$ in $\mathbb{Z}$ and $n \neq 0$) using the DEO families ${\Psi}_{k}^{+}$ and ${\Psi}_{k}^{-}$. In the demonstrations, some close-form formulas are demonstrated such as the decomposition of $f^2$ and $f^3$ with the DEO families ${\Psi}_{k}^{+}$ alone (Lemma) or with  ${\Psi}_{k}^{+}$ and ${\Psi}_{k}^{-}$ (Theorem) if $f$ is chosen in the subset $\mathbf{s}^{-}(\mathbb{R})$. In addition, the theorem shows the existence and uniqueness of the decomposition with the energy operator families, whereas the lemma only shows the existence of the decomposition when using only  $({\Psi}_{k}^{+})_{k \in \mathbb{Z}}$. Note that the lemma and theorem justify the decomposition with the family of energy operators in the case of $f^n$ with $n \in \mathbb{Z}$ and $n> 1$. However, the discussions following the proof of the theorem deal with the cases $n <0$ and $n=1$.
In the following section, the study of the images and kernels helped to simplify some of the formulas decomposing $f^2$ with ${\Psi}_{k}^{+}$ and ${\Psi}_{k}^{-}$. The last part applied this to the energy function and how ${\Psi}_{k}^{+}$ and ${\Psi}_{k}^{-}$ can appear in the development in Taylor-Series of the energy function.
%
\newline This work is an extension of the recently published method to decompose the wave equation with energy operators and a necessary step to extend the method to other linear partial differential equations. 
%
\section*{Acknowledgment}
A special thanks is addressed to Professor Alan McIntosh and Dr. Pierre Portal at the Centre for Mathematics and its Applications at the Australian National University for their inputs and discussions when writing this manuscript. The author also acknowledges the comments from Dr. Ryan Loxton from the Department of Mathematics and Statistics at Curtin University, and from Dr. Herb McQueen from the Research School of Earth Sciences at the Australian National University.
\end{document}